\documentclass[12pt]{article}

\usepackage{authblk}
\usepackage{diagbox}
\usepackage{mathtools}
\usepackage{bbm}
\usepackage{latexsym}
\usepackage{epsfig}
\usepackage{amsmath,amsthm,amssymb,enumerate,hyperref}
\usepackage{comment}
\usepackage[active]{srcltx}
\usepackage{color}

\newcounter{rot}

\def\a{\alpha} \def\b{\beta}

 \def\m{\mu}  
   
 \def\om{\omega}

\newtheorem{theorem}{Theorem}
\newtheorem{lemma}[theorem]{Lemma}

\newcommand{\proofstart}{{\bf Proof\hspace{2em}}}
\newcommand{\proofend}{\hspace*{\fill}\mbox{$\Box$}}
\newcommand{\prooffirst}{{\bf Proof of Theorem \ref{th2}\hspace{2em}}}
\newcommand{\proofsecond}{{\bf Proof of Theorem \ref{th3}\hspace{2em}}}

\newcommand{\brac}[1]{\left(#1\right)}

\newcommand{\bfrac}[2]{\left(\frac{#1}{#2}\right)}

\newcommand{\set}[1]{\left\{#1\right\}}

\def\Pr{\mathbb{P}}

\newcommand{\ignore}[1]{}

\usepackage{tikz}
\usetikzlibrary{arrows}
\usetikzlibrary{decorations}
\usetikzlibrary{shapes.misc}

\def\am{\a_{\min}}

\DeclareMathOperator{\mcp}{mcp}
\DeclarePairedDelimiter{\floor}{\lfloor}{\rfloor}
\begin{document}
\author[1]{Debsoumya Chakraborti\thanks{Supported in part by the Institute for Basic Science (IBS-R029-C1)}}
\author[2]{Mihir Hasabnis}
\affil[1]{\small Discrete Mathematics Group, Institute for Basic Science (IBS), Daejeon,~South~Korea}
\affil[2]{\small Department of Mathematical Sciences, Carnegie Mellon University, Pittsburgh,~USA}
\affil[ ]{\small Email:
\texttt{debsoumya@ibs.re.kr},
\texttt{mhasabni@andrew.cmu.edu} }

\title{The threshold for the full perfect matching color profile in a random coloring of random graphs}
\maketitle

\begin{abstract}{}
Consider a graph $G$ with a coloring of its edge set $E(G)$ from a set $Q = \set{c_1,c_2, \ldots, c_q}$. Let $Q_i$ be the set of all edges colored with $c_i$. Recently, Frieze defined a notion of the perfect matching color profile denoted by $\mcp(G)$, which is the set of vectors $(m_1, m_2, \ldots, m_q)$ such that there exists a perfect matching $M$ in $G$ with $|Q_i \cap M| = m_i$ for all $i$. Let $\a_1, \a_2, \ldots, \a_q$ be positive constants such that $\sum_{i=1}^q \a_i = 1$. Let $G$ be the random bipartite graph $G_{n,n,p}$. Suppose the edges of $G$ are independently colored with color $c_i$ with probability $\alpha_i$. We determine the threshold for the event $\mcp(G) = \set{(m_1, \ldots, m_q) \in [0,n]^q : m_1 + \cdots + m_q = n}$, answering a question posed by Frieze. We further extend our methods to find the threshold for the same event in a randomly colored random graph $G_{n,p}$.  
\end{abstract}

\section{Introduction}
Randomly colored random graphs have been extensively studied in various contexts throughout the last two decades. A few examples include (i) rainbow spanning graphs such as matchings and Hamilton cycles, see e.g., \cite{BF}, \cite{FK}, \cite{FL}, \cite{FNP}, \cite{JW}; (ii) rainbow connection, see e.g., \cite{DFT15}, \cite{HR}, \cite{KKS}, \cite{M}; (iii) pattern colored Hamilton cycles, see e.g., \cite{AF}, \cite{EFK}, \cite{GKM}; (iv) packing problems, see e.g., \cite{FKMS}. Continuing the research in this line, Frieze defined an elegant notion of a color profile in \cite{F} and gave bounds on the matching color profile for randomly colored random bipartite graphs. 

Throughout this paper, we have the following setting: We are given a graph $G$, and positive constants $\a_1, \a_2, \ldots, \a_q$ with $\sum_{i=1}^q \a_i = 1$. Suppose each of the edges of $G$ are independently colored with a random color from the set $Q = \set{c_1, c_2, \ldots, c_q}$ with probability $\Pr\left(c(e) = c_i\right) = \a_i$, where $c(e)$ denotes the color of the edge $e \in E(G)$. Define the color class $Q_i = \set{e \in E(G): c(e) = c_i}$. The perfect matching color profile $\mcp(G)$ is defined to be the set of vectors $(m_1, m_2, \ldots ,m_q)$ such that there exists a perfect matching $M$ in $G$ with $|Q_i \cap M| = m_i$ for all $i$. 

We first consider $G$ to be the random bipartite graph $G_{n,n,p}$. For an event $E_n$, we say that $E_n$ occurs with high probability (in short, w.h.p.) if $\Pr(E_n) \rightarrow 1$ as $n \rightarrow \infty$. Erd\H{o}s and R\'enyi \cite{ER} proved that $G_{n,n,p}$ has a perfect matching w.h.p. when $p=\frac{\log n+\om}{n}$ for any $\om=\om(n)\to\infty$. Moreover, for the same value of $p$, Frieze \cite{F} proved that if the edges of $G = G_{n,n,p}$ are independently colored with $q$ colors with constant probabilities, then most of the elements $(m_1,m_2\ldots m_q) \in [0,n]^q$ such that $\sum_{i=1}^{q} m_i = n$ are present in $\mcp(G)$ w.h.p.

\begin{theorem} [Frieze] \label{th1}
Let $\a_1,\a_2,\ldots,\a_q,\b$ be positive constants such that $\a_1+\a_2+\cdots+\a_q=1$ and $\b<1/q$. Let $G$ be the random bipartite graph $G_{n,n,p}$ where $p=\frac{\log n+\om}{n},\,\om=\om(n)\to\infty$. Suppose that the edges of $G$ are independently colored with colors from $Q=\set{c_1,c_2,\ldots,c_q}$ where $\Pr(c(e)=c_i)=\a_i$ for $e\in E(G),i\in[q]$. Let $m_1,m_2,\ldots,m_q$ satisfy: (i) $m_1+\cdots+m_q=n$ and (ii) $m_i\geq \b n,i\in [q]$. Then w.h.p., there exists a perfect matching $M$ in which exactly $m_i$ edges are colored with $c_i,i=1,2,\dots,q$.
\end{theorem} 

It is not hard to check that w.h.p. $(n ,0 , \ldots, 0) \notin \mcp(G)$, in view of the fact that the bipartite graph induced by the first color is distributed as $G_{n,n,\a_1p}$ and has isolated vertices w.h.p. Frieze posed the natural problem of determining the threshold for $\mcp(G) = \set{(m_1, \ldots, m_q) \in [0,n]^q : m_1 + \cdots + m_q = n}$. In this paper, we determine that threshold. 

\begin{theorem}\label{th2}
Let $\a_1,\a_2,\ldots,\a_q$ be positive constants such that $\a_1+\a_2+\cdots+\a_q=1$. Let
\[
\am=\min\set{\a_i:i\in [q]}.
\] 
Let $G$ be the random bipartite graph $G_{n,n,p}$ where $p=\frac{\log n+\om}{\am n},\,\om=\om(n)\to\infty$. Suppose that the edges of $G$ are independently colored with colors from $C=\set{c_1,c_2,\ldots,c_q}$ where $\Pr(c(e)=c_i)=\a_i$ for $e\in E(G),i\in[q]$. Then, w.h.p. for each $m_1,m_2,\ldots,m_q$ satisfying $m_1+\cdots+m_q=n$, there exists a perfect matching $M$ in which exactly $m_i$ edges are colored with $c_i,i=1,2,\dots,q$. In other words, $\mcp(G) = \set{(m_1, \ldots, m_q) \in [0,n]^q : m_1 + \cdots + m_q = n}$.
\end{theorem} 

Let us first determine the lower bound on the threshold. Assume that $\am = \a_i$. To prove the lower bound, note that it is enough to show that the same threshold holds even for the event that $G$ contains a perfect matching in color $c_i$. To see this, remember that the bipartite graph induced by the color $c_i$ is distributed as $G_{n,n,\a_ip}$. The claim now follows from the known thresholds of the random bipartite graph to have a perfect matching, see e.g., Theorem 6.1 of \cite{FK0}. The general strategy to prove the upper bound on the threshold in Theorem \ref{th2} is to do the following modification iteratively. For each $i \neq j$, if $G$ contains a perfect matching $M$ using $m_i \ge \frac{n}{q}$ edges with color $c_i$, then we can find a perfect matching $M'$ consisting of one fewer edge of color $c_i$ and one more edge of color $c_j$. Frieze \cite{F} also suggested studying the same problem for the random graph $G_{n,p}$. A simple extension of our techniques establishes the threshold for $G_{n,p}$ as well.

\begin{theorem}\label{th3}
Let $\a_1,\a_2,\ldots,\a_q$ be positive constants such that $\a_1+\a_2+\cdots+\a_q=1$. Let
\[
\am=\min\set{\a_i:i\in [q]}.
\] 
Let $G$ be the random graph $G_{n,p}$ where $p=\frac{\log n+\om}{\am n},\,\om=\om(n)\to\infty$. Suppose that the edges of $G$ are independently colored with colors from $C=\set{c_1,c_2,\ldots,c_q}$ where $\Pr(c(e)=c_i)=\a_i$ for $e\in E(G),i\in[q]$. Then, w.h.p. for each $m_1,m_2,\ldots,m_q$ satisfying $m_1+\cdots+m_q=\floor{\frac{n}{2}}$, there exists a perfect matching $M$ in which exactly $m_i$ edges are colored with $c_i,i=1,2,\dots,q$. In other words, $\mcp(G) = \set{(m_1, \ldots, m_q) \in [0,n]^q : m_1 + \cdots + m_q = \floor{\frac{n}{2}}}$.
\end{theorem} 

Similar to Theorem \ref{th2}, the lower bound on the threshold for Theorem \ref{th3} follows from the known thresholds of the random graph to have a perfect matching (see, e.g., Theorem 6.2 of \cite{FK0}).

This short note is organized as follows. The next section is devoted to stating a few simple structural lemmas about random bipartite graphs and random graphs. Section 3 contains the proof of Theorem \ref{th2} and Theorem \ref{th3}. Finally, we finish with a few concluding remarks.

\section{Structural lemmas}
Let $\a_i, 1 \le i \le q$, and $\am$ be as in Theorems \ref{th2} and \ref{th3}. Throughout this section, the graph $G$ will be either the random bipartite graph $G_{n,n,p}$ or the random graph $G_{n,p}$, where the probability $p = \frac{\log n + \omega}{\am n}$, for some $\omega = \omega(n) \rightarrow \infty$. The edges of $G$ are randomly colored as in Theorems \ref{th2} and \ref{th3}.

\begin{lemma}\label{lem1}
Let $G$ be the random bipartite graph $G_{n,n,p}$ with the vertex bipartition $A \cup B$. Suppose that the edges of $G$ are independently colored with colors from $C = \{c_1, c_2, \ldots, c_q\}$ where each edge is colored with $c_i$ by probability $\a_i$. Then, w.h.p. for each $i \in [q]$, and any $X \subseteq A, Y \subseteq B$ with $|X|, |Y| \geq \frac{n}{4q}$, there is an edge with color $c_i$ between $X$ and $Y$ in $G$.
\end{lemma} 

\proofstart 
Note that it is enough to prove this lemma with $|X| = |Y| = \frac{n}{4q}$. Now by a simple union bound, we have the following:
\begin{align*}
\mathbb{P}(\exists X, Y \text{	s.t. condition is not satisfied}) &\le \binom{n}{n/4q}^2 \sum_{i=1}^q \brac{1 - p\a_i}^{\frac{n^2}{16q^2}} \\
&\le q \bfrac{ne}{n/4q}^{n/2q} \brac{1 - \frac{\log n}{n}}^{\frac{n^2}{16q^2}} \\
&\le q \brac{(4eq)^{1/2q} \cdot e^{-\frac{\log n}{16q^2}}}^n \\
&= o(1).
\end{align*}
\proofend

\begin{lemma} \label{lem2}
Let $G$ be the random graph $G_{n,p}$. Suppose that the edges of $G$ are colored in the exact same way as in Lemma \ref{lem1}. Then, w.h.p. for each $i \in [q]$, and any disjoint sets $X, Y \subseteq V(G)$ with $|X|, |Y| \geq \frac{n}{8q}$, there is an edge with color $c_i$ between $X$ and $Y$ in $G$.
\end{lemma}

\proofstart
This follows very similarly to the proof of Lemma \ref{lem1}.
\proofend

\begin{lemma} \label{lem3}
Let $G$ be the random bipartite graph $G_{n,n,p}$ or the random graph $G_{n,p}$. Then, w.h.p. for each $i \in [q]$, the graph $G$ contains a perfect matching in color $c_i$.
\end{lemma}

\proofstart
This is an easy consequence of Theorems 6.1 and 6.2 of \cite{FK0}.
\proofend

\section{Proof of the main results}
\prooffirst
Suppose that we are given a bipartite graph $G$ for which the high probability properties (Lemmas \ref{lem1} and \ref{lem3}) of the random bipartite graph $G_{n,n,p}$ mentioned in the last section hold. The proof mainly consists of showing that the following can be done. For each $i \neq j$, if $G$ contains a perfect matching $M$ with at least $\frac{n}{q}$ edges with color $c_i$, then $G$ contains a perfect matching with the same color profile as $M$ but with one fewer edge of color $c_i$ and one more edge of color $c_j$. We next show how we can iteratively apply this modification to obtain a perfect matching with any given color profile.

Fix $(m_1, m_2, \ldots, m_q) \in [0,n]^q$ such that $\sum_{i=1}^q m_i = n$. Our goal is to show that $G$ has a perfect matching $M$ such that $|M \cap Q_i| = m_i$ for all $i$. Without loss of generality we can assume that $m_1 = \max\set{m_i:i\in [q]}$. This implies that $m_1 \ge \frac{n}{q}$. By Lemma \ref{lem3}, we know that there is a perfect matching in the subgraph induced by color $c_1$ in $G$. We proceed in the following way: starting with a perfect matching with color profile $(n, 0, \ldots, 0)$, for any fixed color $c_j$ with $j \neq 1$ we show the existence of a perfect matching with one fewer edge in color $c_1$ and one more edge in color $c_j$. We keep doing this process until we get a matching with $m_i$ edges with color $c_i$ for all $i$. Note that we need $n - m_1$ steps to reach a matching with the color profile $(m_1, m_2, \ldots, m_q)$, because in every step, we find a matching with one fewer edge in color $c_1$. So, it is enough to show that for any perfect matching $M$ in $G$ with $|M \cap Q_i| = \m_i$ for each $i \in [q]$ and $\m_1 \ge \frac{n}{q}$, there is a matching $M'$ in $G$ with $|M' \cap Q_1| = \m_1 - 1$, $|M' \cap Q_2| = \m_2 + 1$ and $|M' \cap Q_i| = \m_i$ for all other $i$. 

We show the above statement by finding an appropriate alternating cycle. More precisely, we find a cycle $C$ with vertex sequence $(x_1 \in A, y_1 \in B, x_2 \in A, y_2 \in B, \ldots, x_\ell \in A, y_\ell \in B, x_1)$ such that (i) $(x_i, y_i) \not \in M$, (ii) $(y_i, x_{i+1}) \in M$, (iii) $(x_1, y_1) \in Q_2$, and (iv) $E(C) \setminus \set{(x_1, y_1)} \subseteq Q_1$. For the convenience of writing the proof, we introduce some notation. Label vertices so that the edges $v_i^+ v_i^-$, $i \in [\frac{n}{q}]$, with $v_i^+ \in A$ and $v_i^- \in B$ are distinct edges with color $c_1$ in $M$. Create a directed graph $D$ on vertex set $\set{v_1, \ldots, v_{n/q}}$, where there is a directed edge $v_i v_j$ in $D$ if there is an edge with color $c_1$ between $v_i^-$ and $v_j^+$ in $G$.

Note that if there is an edge with color $c_2$ between $v_i^+$ and $v_j^-$ in $G$ and a directed path from $v_i$ to $v_j$ in $D$, then this gives exactly the alternating cycle $C$ which we discussed in the last paragraph. Furthermore, by using Lemma \ref{lem1}, we have the following property in $D$. 
\begin{enumerate}
\item For each $X, Y \subseteq V(D)$ with $|X|, |Y| \ge \frac{n}{4q}$, there is an edge from $X$ to $Y$ in $D$. 
\end{enumerate}

For each $v \in V(D)$, let $B^+(v)$ be the set of vertices reachable by a directed path from $v$ in $D$ (including $v$), and let $B^-(v)$ be the set of vertices in $V(D)$ from which you can reach $v$ in $D$ with a directed path (including $v$). Let $V_1 = \set{v \in V(D) : |B^+(v)| \le \frac{n}{4q}}$ and $V_2 = \set{v \in V(D) : |B^-(v)| \le \frac{n}{4q}}$.

Now, claim that $|V_1| \le \frac{n}{4q}$. If not, then we can pick a minimal set $V_1' \subseteq V_1$ such that $|\cup_{v \in V_1'} B^+(v)| \ge \frac{n}{4q}$, and note that $|\cup_{v \in V_1'} B^+(v)| \le \frac{2n}{4q}$. There are no edges from $\cup_{v \in V_1'} B^+(v)$ into $V(D) \setminus \left(\cup_{v \in V_1'} B^+(v)\right)$, and the latter set has size at least $|D| - \frac{2n}{4q} \ge \frac{n}{4q}$, contradicting the property (1). Therefore, $|V_1| \le \frac{n}{4q}$. Similarly, $|V_2| \le \frac{n}{4q}$. Thus, $|V(D) \setminus (V_1 \cup V_2)| \ge \frac{n}{2q}$.    

Now, by Lemma \ref{lem1} there is an edge in $G$ with color $c_2$ between $\set{v_i^+ : v_i \in V(D) \setminus (V_1 \cup V_2)}$ and $\set{v_i^- : v_i \in V(D) \setminus (V_1 \cup V_2)}$. Say this is the edge $v_i^+ v_j^-$ and note that $i \neq j$. As $v_i, v_j \in V(D) \setminus (V_1 \cup V_2)$, we have that $|B^+(v_i)|, |B^-(v_j)| \ge \frac{n}{4q}$. Thus, there is an edge from $B^+(v_i)$ into $B^-(v_j)$ in $D$ by (1), and therefore there is a directed path from $v_i$ to $v_j$ in $D$. This finishes the proof of Theorem \ref{th2}.
\proofend

\proofsecond
The proof of Theorem \ref{th2} extends straightforwardly to a proof of Theorem \ref{th3}. By Lemma \ref{lem3}, we know that $G = G_{n,p}$ has a perfect matching in each color. Now, if a color profile $(m_1, \ldots, m_q)$ is required (say $m_1$ is the largest of these), then start with a perfect matching in color $c_1$, and split $V(G)$ into $A$ and $B$ arbitrarily so that $M$ is a matching between $A$ and $B$. The same arguments as in the proof of Theorem \ref{th2} can now be used due to Lemma \ref{lem2}, which is the replacement of Lemma \ref{lem1} we used before. More precisely, to modify a perfect matching $M$ to another matching $M'$ with the same color profile but one fewer edge of color $c_1$ and one more edge of color $c_j$, we choose an arbitrary bipartition $V(G) = A \cup B$ with $M$ being a matching between $A$ and $B$, and then implement the exact same argument as before. 
\proofend

\section*{Concluding remarks}
In this short note, we consider the random bipartite graph $G = G_{n,n,p}$ and the random graph $G_{n,p}$, and determine the threshold on the parameter $p$ for the event that $G$ contains perfect matchings of all color profiles, i.e., $\mcp(G) = \set{(m_1, \ldots, m_q) \in [0,n]^q : m_1 + \cdots + m_q = n}$. Some interesting directions of future research would be to determine $\mcp(G)$ for Hamilton cycles, spanning trees etc. or to consider deterministic host graphs (e.g., Dirac graphs) instead of random graphs.

\section*{Acknowledgements}
The authors are grateful to Alan Frieze for introducing this problem to them. The authors are indebted to the anonymous reviewer for simplifying the proof of Theorem \ref{th2} in the original version of this paper.


\begin{thebibliography}{99}
\bibitem{AF} M. Anastos and A.M. Frieze, Pattern Colored Hamilton Cycles in Random Graphs, {\em SIAM Journal on Discrete Mathematics} 33 (2019) 528-545.
\bibitem{BF} D. Bal and A.M. Frieze, Rainbow Matchings and Hamilton Cycles in Random Graphs, {\em Random Structures and Algorithms} 48 (2016) 503-523.
\bibitem{ER} P. Erd\H{o}s and A. R\'enyi, On random matrices, {\em Publ. Math. Inst. Hungar. Acad. Sci.} 8 (1964) 455-461.
\bibitem{DFT15} A. Dudek, A.M. Frieze and C.E. Tsourakakis, Rainbow connection of random regular graphs, {\em SIAM Journal on Discrete Mathematics} 29 (2015) 2255-2266.
\bibitem{EFK} L. Espig, A.M. Frieze and M. Krivelevich, Elegantly colored paths and cycles in edge colored random graphs, {\em SIAM Journal on Discrete Mathematics} 32 (2018) 1585-1618.
\bibitem{F} A.M. Frieze, A note on randomly colored matchings in random graphs, Discrete Mathematics and Applications, Springer Optimization and Its Applications, 165 (2020).
\bibitem{FK0} A. Frieze and M. Karo\'nski, Introduction to Random Graphs, {\em Cambridge University Press} 2015.
\bibitem{FK} A. Ferber and M. Krivelevich, Rainbow Hamilton cycles in random graphs and hypergraphs, Recent trends in combinatorics, IMA Volumes in Mathematics and its applications, A. Beveridge, J. R. Griggs, L. Hogben, G. Musiker and P. Tetali, Eds., Springer 2016, 167-189.
\bibitem{FKMS} A. Ferber, G. Kronenberg, F. Mousset, and C. Shikhelman, Packing a randomly
edge-colored random graph with rainbow k-outs, {\em arXiv preprint arXiv:1410.1803} (2014).
\bibitem{FL} A.M. Frieze and P. Loh, Rainbow Hamilton cycles in random graphs, {\em Random Structures and Algorithms} 44 (2014) 328-354. 
\bibitem{FNP} A. Ferber, R. Nenadov, and U. Peter, Universality of random graphs and rainbow embedding. {\em Random Structures and Algorithms} 48(3) (2016) 546-564.
\bibitem{GKM} L. Gishboliner, M. Krivelevich, and P. Michaeli, Colour-biased Hamilton cycles in random graphs, {\em arXiv preprint arXiv:2007.12111} (2020).
\bibitem{HR} A. Heckel and O. Riordan, The hitting time of rainbow connection number two, {\em Electronic Journal on Combinatorics} 19 (2012).
\bibitem{JW} S. Janson and N. Wormald, Rainbow Hamilton cycles in random regular graphs,
{\em Random Structures Algorithms} 30 (2007) 35-49.
\bibitem{KKS} N. Kamcev, M. Krivelevich and B. Sudakov, Some remarks on rainbow connectivity, {\em Journal of Graph Theory} 83 (2016) 372-383. 
\bibitem{M} M. Molloy, The rainbow connection number for random 3-regular graphs,  {\em Electronic Journal of Combinatorics} 24 (2017).
\end{thebibliography}
\end{document}